\documentclass{amsart}
\usepackage[foot]{amsaddr}
\usepackage[utf8]{inputenc}
\usepackage[T1]{fontenc}
\usepackage{amssymb}
\usepackage{amsmath}
\usepackage{mathtools}
\usepackage{geometry}
\makeatletter
\DeclareRobustCommand*\cal{\@fontswitch\relax\mathcal}

\usepackage{amscd,amsmath,amssymb,amsthm,amsfonts,epsfig,graphics}
\newtheorem{thrm}{Theorem}[section]
\newtheorem{prpstn}{Proposition}
\newtheorem{lmm}{Lemma}

\newtheorem{dfntn}{Definition}
\newtheorem{rmrk}{Remark}

\title{Hold-out estimates of prediction models for Markov processes}

\thanks{This work was partly involved in Project ”EcoDep” PSI-AAP2020-0000000013}

\author{Remy Garnier$^1$}
\address[1]{CDiscount\\
  120-126 Quai de Bacalan, 33300 Bordeaux, France}

\author{Rapha\"el Langhendries$^2$}
\author{Joseph Rynkiewicz$^2$}
\address[2]{Universit\'e Paris I - SAMM\\
90 rue de tolbiac, Paris - France}

\begin{document}
\maketitle

\begin{abstract}
We consider the selection of prediction models for Markovian time series. For this purpose, we study the theoretical properties of the hold-out method. In the econometrics literature, the hold-out method is called ``out-of-sample'' and is the main method to select a suitable time series model.  This method consists of estimating models on a learning set and picking up the model with minimal empirical error on a validation set of future observations. Hold-out estimates are well studied in the independent case, but, as far as we know, this is not the case when the validation set is not independent of the learning set. In this paper, assuming uniform ergodicity of the Markov chain, we state generalization bounds and oracle inequalities for such method;  in particular, we show that the ``out-of-sample'' selection method is adaptative to noise condition.
\end{abstract}



\section{Introduction}
Many models in time series involve the one-step prediction of the next value knowing past values, and hold-out or out-of-sample (OOS) method is probably the most commonly used model selection method in practice. This method consists in splitting the sample of size $n+m$ in two parts: a training set of length $n$ and a validation set of size $m$. The training set is used to derive a finite collection of candidate prediction functions and we chose the function with the best performance on the validation set. It turns out to look like picking a prediction function from a finite collection; however, we must be careful with the dependence between the learning set and the validation set for the Markov case. For time series, splitting the data into a training subset and a validation subset of future observations is an option implemented in most statistical or machine learning software. This standard evaluation procedure works very well in practice (see, for example, Tashman \cite{Tashman} for OOS to assess the model's accuracy). In the machine learning community, Cerqueira et al. \cite{Cerqueira} compare the performances between OOS and other methods empirically; they found that OOS produces the most accurate estimates for real time series. The authors think that the main reason for the performance of OOS method is the preservation of the temporal order of the observations. For the independent and identical distributed (i.i.d.) case, the hold-out theoretical properties are well known, and, for example, in the classification case, it adapts to the noise conditions (see Blanchard and Massart \cite{Blanchard}).
However, as far as we know, there are few theoretical results for the OOS method for dependent data. Some studies assess the asymptotical performance of related methods like cross-validation (see Arlot and Celisse \cite{Arlot}) in the context of regression (Burman and Nolan \cite{Burman}). But, for model selection, in the case of dependent observations, cross-validation is known to be severely affected by dependence (see Chu and Marron \cite{Chu}). Other authors investigate methods of model selection for dependent data, which focus more on penalization (see Alquier and Wintenberger \cite{Alquier}) or on complexity measure such that Rademacher complexity (Mohri and Kuznetsov \cite{MohriRad}) or stability bounds (Mohri and Rostamizedeh \cite{MohriStab}). Empirical risk minimization has also been studied in the framework of uniformly ergodic Markov chains (see Bin et al. \cite{Bin}). To be exhaustive, we can also cite some studies on the asymptotic convergence rate for estimating finite Markov chain transition matrices (see Falahatgar et al. \cite{Falahatgar}, Hao et al. \cite{Hao}, and \cite{Han}). Note that such results are difficult to apply to massive models like Deep networks (Zhang et al. \cite{Zhang}). Indeed, these authors show, through extensive systematic experiments, that these traditional approaches fail to explain why large neural networks generalize well in practice. Eventually, OOS is still the standard method for model selection in time series, especially for Deep learning models.

If the data are drawn from a process indexed in time order, the validation set is no longer independent from the learning set, and the classical i.i.d. theory of hold-out does not hold anymore. This paper aims to provide generalization bounds and oracle inequalities for the selected prediction model by the OOS method in a Markovian framework.
Our paper is organized as follows: in the next section, we will present the observations, models, notations, and concentration inequalities for uniformly ergodic Markov chain. In the third section, we give first exponential inequalities, generalization bounds, and oracles inequalities for the model selected with OOS method. In this section, we use only the boundedness property of the loss function, and under additional assumptions, we improve these bounds in the fourth section. Finally, in the fifth section, we refine these bounds under noise conditions and show our main result: the OOS method is still adaptative to noise conditions for uniformly ergodic Markov chain. We postpone long proofs in the Appendix.

\section{The model}
\subsection{Assumptions and definitions for the observed process}
We consider $(Y_t)_{t\in\mathbb Z}$, an $\cal Y$-valued, $k$-order Markov chain, where ${\cal Y}$ is a Polish state space. We assume that $(Y_t)_{t\in\mathbb Z}$ is  time homogenous, stationary, and uniformly ergodic.  Let $p$ be fixed a integer with $p\geq k$, then the markovization $(X_t)_{t\in\mathbb Z}$ of $(Y_t)_{t\in\mathbb Z}$ with:
\begin{equation}\label{defXt}
X_t\coloneqq (Y_{t},\cdots,Y_{t-p})^T,
\end{equation}
will be a time homogenous, stationary, uniformly ergodic Markov chain of order $1$.
We are going to state some definitions from the theory of general state space Markov chains, based on Gareth and Rosenthal \cite{Gareth}.
\begin{dfntn}\label{BasicDefMarkov}
Let us denote $K$ its transition kernel. $K(x,.)$ is the distribution of $X_{n+1}$ conditioned on $X_n=x$. Denote by $Q$ its stationary law, i.e. the probability distribution such that
\begin{equation*}
\int_{x\in\cal X}Q(dx)K(x,dz)=Q(dz).
\end{equation*}
We define the total variational distance of two distributions $P$ and $Q$ defined on the same state space $(\cal X,\cal A)$ as
\begin{equation}\label{totalvariationsup}
  d_{TV}(P,Q)\coloneqq \sup_{A\in\cal A}|P(A)-Q(A)|.
\end{equation}

Since $(Y_t)_{t\in\mathbb Z}$ is uniformly ergodic, $(X_t)_{t\in\mathbb Z}$ will be uniformly ergodic, and constants $C>0$ and $0< \rho<1$ exist such that:
\begin{equation} \label{Uniformergodicity}
  \sup_{x\in\cal X}d_{TV}(K^n(x,.),Q)\leq C\rho^n.
\end{equation}
\end{dfntn}\label{defMarkov}
\subsection{Estimated functions and loss functions}
We want to estimate a one-step prediction model:  $g^*(Y_{t-1},\cdots,Y_{t-p})$, which minimizes the expectation of a real bounded loss function $L$.

According to equation \eqref{defXt}, the observations are a realization of
\(\left(X_1,\cdots,X_n,X_{n+1},\cdots,X_{n+m}\right)\).
the variables $\left(X_1,\cdots,X_n\right)$ constitute the learning set, and the variables $(X_{n+1},\cdots,X_{n+m})$ the validation set. Let us introduce some definitions:
\begin{dfntn}\label{Bestfunction} We will introduce the notion of prediction and loss functions.
  \begin{itemize}
   \item A  measurable function $g$ from ${\cal Y}^{p}$ into $\cal Y$ will be called a prediction function. 
     \item Let $L$ be a real, positive, measurable, bounded function defined on ${\cal Y}^2$. Without loss of generality, we can always rescale the function $L$ such that $|L(y,y')|\leq 1$. $L$ will be called a loss function.
  \item With a slight abuse of notation, we will denote:
\begin{equation}\label{LgX}
  L(g(X_t))\coloneqq L\left(g(Y_{t-1},\cdots,Y_{t-p}),Y_{t}\right).
\end{equation}
  \item  Let $X$ be a random vector with the stationary law $Q$ of the Markov chain $(X_t)_{t\in\mathbb Z}$. The expected loss of a measurable function $g$ will be:
\begin{equation}\label{defL}
  {\mathbb L}(g)={\mathbb E}_Q\left(L(g(X))\right).
\end{equation}
  \item Let $\cal F$ be the set of measurable functions from ${\cal Y}^{p}$ into $\cal Y$. The best prediction function, $g^*$, is the function that minimizes the loss function applied to the one-step prediction under the stationary law:
\begin{equation}
g^*=\arg\min_{g\in \cal F}{\mathbb E}_QL\left(g(X_t)\right).
\end{equation}
\end{itemize}
\end{dfntn}
For example, if $\cal Y$ is finite, and  $L$ is the misclassification loss function: $L(y,y')={\bf 1}_{y\neq y'}$, $g^*$ will the best prediction of the following state $Y_{t}$ knowing $Y_{t-1},\cdots,Y_{t-p}$.

To estimate $g^*$, we seek prediction function $\hat g$ among a set of possible functions ${\cal G}\subset{\cal F}$ by minimizing an empirical loss function on the learning set $(X_{1},\cdots,X_n)$:
\begin{dfntn}\label{Estimatedfunction} 
  Let $\tilde L$ be a loss function.
  \begin{itemize}
    \item The empirical estimation of $g^*$, among a set of possible functions ${\cal G}$, will be defined as:
  \begin{equation}\label{argming}
  \hat g\left(X_1,\cdots,X_n\right)=\arg\min_{g\in{\cal G}}\sum_{t=1}^{n}{\tilde L}\left(g(X_{t})\right).
  \end{equation}
  \item The estimated function $\hat g$ depends on $X_1,\cdots,X_n$ and is a random function. If we observe the realization $x_1,\cdots,x_n$ of $X_1,\cdots,X_n$, we will observe a realization $\hat g(x_1,\cdots,x_n)$ of $\hat g(X_1,\cdots,X_n)$. For convenience, since our results are valid for any realization $\hat g(x_1,\cdots,x_n)$ of $\hat g(X_1,\cdots,X_n)$, we will denote
\begin{equation}
  {\hat g}_1^n\coloneqq \hat g(x_1,\cdots,x_n).
\end{equation}
\end{itemize}
\end{dfntn}

 The set of possible function ${\cal G}$ depends on user choice. In pratice, we have to chose between several sets $\left\{{\cal G}_k\right\}_{1\leq k\leq N}$, where $N$ is a finite integer. Each set ${\cal G}_k$ defines a empirical minimizer:
\begin{equation}
   ({\hat g}_1^n)_k=\arg\min_{g\in{\cal G}_k}\sum_{t=1}^{n}{\tilde L}\left(g(x_{t})\right).
\end{equation}
Our goal is to chose the best prediction function among $\left(({\hat g}_1^n)_k\right)_{1\leq k\leq N}$.

Note that the learning loss $\tilde L$ need not be equal to $L$; it can also be a proxy function easier to optimize.
If the set of possible functions ${\cal G}_k$ is large, the function $({\hat g}_1^n)_k$ may have poor performances on future observations of the process $(X_t)_{t\in\mathbb Z}$, even if the empirical learning loss is small. We say that $({\hat g}_1^n)_k$ overfits the learning set $(X_{1},\cdots,X_n)$.  Hence, it is wise to assess the performance of the estimated model on an out-of-sample set $X_{n+1},\cdots,X_{n+m}$:
\begin{dfntn}\label{theoloss}
Let $g$ be any measurable function from ${\cal Y}^p$ into $\cal Y$. The empirical, out-of-sample loss of $g$ will be: 
	\begin{equation}\label{empiricalL}
	\hat L_{m}(g)=\frac{1}{m}\sum_{k=n+1}^{n+m}L(g(X_k)).
	\end{equation}
\end{dfntn}

Note that, by the law of large number, the empirical loss ${\hat L}_m(g)$ converges, almost surely, towards the theoretical loss ${\mathbb L}(g)$:
        \[
{\mathbb L}(g)\stackrel{a.s.}{=}\lim_{m\rightarrow\infty}{\hat L}_m(g).
\]

\begin{rmrk}
For an integer $b\geq 0$, let us write ${\mathbb L}_b({\hat g}_1^n)$, the expected loss of any minimizer ${\hat g}_1^n$ for a vector $X_{n+1+b}$ of future observations conditionally to the realization $(X_1=x_1,\cdots,X_n=x_n)$:
\begin{equation}
{\mathbb L}_b({\hat g}_1^n)={\mathbb E}\left(L({\hat g}_1^n(X_{n+b+1}))\left|X_1=x_1,\cdots,X_n=x_n\right.\right).
\end{equation}
A direct application of the uniform ergodicity of the Markov chain shows that, for a finite $b$, we can approximate the expectation of the loss ${\mathbb L}_b$ by the theoretical loss ${\mathbb L}$. Hence, with the notations of definition \ref{BasicDefMarkov}, for any realization $x_1,\cdots,x_n$ of $X_1,\cdots,X_n$:
\begin{equation}\label{couplage}
  \left|{\mathbb L}_b({\hat g}_1^n)-{\mathbb L}({\hat g}_1^n)\right|\leq C\rho^b.
\end{equation}
\end{rmrk}
\subsection{Exponential inequalities}
To study the link between the empirical loss \eqref{empiricalL} and the theoretical loss \eqref{defL}, we need uniform inequalities between the empirical mean and the expected mean (as in Lugosi \cite{LugosiIB}).  This section aims to give such inequalities. 

First, let us introduce the notion of mixing time, which allows us to evaluate the speed of the convergence of a Markov Chain to its stationary distribution. 

\begin{dfntn}\label{Mix_time}
  Let $\left(X_t\right)_{t\in\mathbb Z}$ be a time homogeneous, uniformly ergodic, Markov chain. Let the total variation distance be defined by equation \eqref{totalvariationsup}. The mixing time $t_{mix}$ is defined by:
	\[d(t)\coloneqq \underset{x\in \cal X}{\sup} \| K^t(x, \cdot) - Q \|_{TV},\ t_{mix}(\varepsilon) = \min \{t: d(t) \leq \varepsilon \},\mbox{ and } t_{mix} = t_{mix}(\frac{1}{4}).\]
\end{dfntn}
The fact that $t_{mix}$ is finite is equivalent to the uniform ergodicity of the chain (see Roberts and Rosenthal \cite{Roberts}). 
If we introduce an integer gap $b$ between the learning and validation sample as a technical tool, a straightforward adaptation of Corollary 2.10 and equation (3.27) of Paulin \cite{Paulin} yields the following proposition:
\begin{prpstn}\label{HoeffdingPaulin}
  Let $C\geq 1$ and $0\leq \rho < 1$ be the positive constants in equation  \eqref{Uniformergodicity}. For any realization $x_1,\cdots,x_n$ of $X_1,\cdots,X_n$, real number $0\leq \varepsilon\leq 1$, integers $0\leq b<m$, with the notations of definitions  \ref{theoloss}, \ref{Estimatedfunction}, and \ref{Mix_time}:
  \begin{equation}\label{Hoeffdingprel}
      {\mathbb P}\left(\pm\left(\frac{1}{m-b}\sum_{k=n+1+b}^{n+m}L({\hat g}_1^n(X_k))-{\mathbb L}({\hat g}_1^n)\right)>\varepsilon\right)
      \leq \exp\left(-2\frac{(m-b)\varepsilon^2}{9t_{mix}}\right)+C\rho^{b}.
  \end{equation}
  Moreover, according to equation (3.30) of Paulin \cite{Paulin}, we have $C\rho^{b}\leq 2\exp\left(-\frac{b\ln(2)}{t_{mix}}\right)$, so
    \begin{equation}\label{Hoeffding}
      {\mathbb P}\left(\pm\left(\frac{1}{m-b}\sum_{k=n+1+b}^{n+m}L({\hat g}_1^n(X_k))-{\mathbb L}({\hat g}_1^n)\right)>\varepsilon\right)\leq
       \exp\left(-2\frac{(m-b)\varepsilon^2}{9t_{mix}}\right)+2\exp\left(-\frac{b\ln(2)}{t_{mix}}\right).
  \end{equation}
\end{prpstn}

Now, we introduce the notion of pseudo spectral gap briefly; see Paulin \cite{Paulin} for a detailed presentation.
\begin{dfntn}\label{Spectral_gap}
  For a Markov chain with transition kernel $K(x,dz)$ and stationary distribution $Q$, we define the spectrum of the chain as
  \[
    S_2\coloneqq \left\{\lambda\in {\mathbb C}\backslash 0:(\lambda {\bf I}-K)^{-1}\mbox{ does not exist as a bounded linear operator on }
                                                                          L^2(Q)\right\}.
   \]
   We also define the time reversal of $K$ as the Markov kernel
   \[
   K^*(x,dz)\coloneqq \frac{K(z,dx)}{Q(dx)}Q(dz).
   \]
   Then, the linear operator $K^*$ is the adjoint of the linear operator $K$ on $L^2(Q)$. We define a new quantity, called the pseudo spectral gap of $K$, as
   \begin{equation*}
     \gamma_{ps}\coloneqq \max_{k\geq 1}\left\{\gamma\left((K^*)^kK^k\right)/k\right\},
   \end{equation*}
   where $\gamma\left((K^*)^kK^k\right)$ denotes the spectral gap of the self-adjoint operator $(K^*)^kK^k$.
\end{dfntn}
Now, a straightforward adaptation of Theorem 3.4 of Paulin \cite{Paulin} gives:

\begin{prpstn}\label{PreBernsteinPaulin}
  Let $(X_t)_{t\in\mathbb N}$ be a stationary Markov chain with spectral gap $\gamma_{ps}$, and a real number $0\leq \varepsilon\leq 1$. Let $f\in {L^2(Q)}$ with, for every $x$, $|f(x)-E_Q(f)|\leq B$. Let $V_f=Var_Q(f)$, and $S=\sum_{i=1}^nf(X_i)$, then:
    \begin{equation}\label{PreBernstein}
      {\mathbb P}\left(\pm\left(S-E_Q(S)\right)>n\varepsilon\right)\leq
      \exp\left(-\frac{n^2\varepsilon^2\gamma_{ps}}{8(n+1/\gamma_{ps})V_f+20n\varepsilon B}\right).
  \end{equation}
\end{prpstn}
From this proposition, we deduce a lemma, proven in the Appendix, that will be used in the last section:
\begin{lmm}\label{lemmemoinsdelta}
  With the same assumptions than the previous proposition \ref{PreBernsteinPaulin}, for any $0<\delta<1$:
  \begin{align}\label{borne1moinsdelta}
      & P\left(\pm\left(E_Q(S)-S\right)\leq \sqrt{\frac{8\left(\gamma_{ps}+1\right)}{\gamma^2_{ps}}nV_f\log\left(\frac{1}{\delta}\right)}+\frac{20}{\gamma_{ps}}B\log\left(\frac{1}{\delta}\right)\right)\geq \nonumber\\
      & P\left(\pm\left(E_Q(S)-S\right)\leq \sqrt{\frac{8}{\gamma_{ps}}(n+1/\gamma_{ps})V_f\log\left(\frac{1}{\delta}\right)}+\frac{20}{\gamma_{ps}}B\log\left(\frac{1}{\delta}\right)\right)\geq 1-\delta
    \end{align}
\end{lmm}
Finally,  proposition \ref{PreBernsteinPaulin} and equation (3.27) of Paulin \cite{Paulin} yield the following proposition:
\begin{prpstn}\label{BernsteinPaulin}
  Let $(X_t)_{t\in\mathbb Z}$ be a stationary Markov chain with pseudo spectral gap $\gamma_{ps}$. For any function $g$,  let us denote $V_g=Var(L(g(X_t))$ the variance of the loss function computed with the stationary law.  Then, for any realization $x_1,\cdots,x_n$ of $X_1,\cdots,X_n$, real number $0\leq \varepsilon\leq 1$, integers $0\leq b<m$, with the notations of definitions  \ref{theoloss} and  \ref{Estimatedfunction}, since  $C\rho^{b}\leq 2\exp\left(-\frac{b\ln(2)}{t_{mix}}\right)$:
    \begin{align}\label{Bernstein}
     & {\mathbb P}\left(\pm\left(\frac{1}{m-b}\sum_{k=n+b+1}^{n+m}L({\hat g}_1^n(X_k))-{\mathbb L}({\hat g}_1^n)\right)>\varepsilon\right)\leq \nonumber \\
     & \exp\left(-\frac{(m-b)^2\varepsilon^2\gamma_{ps}}{8((m-b)+\frac{1}{\gamma_{ps}})V_{{\hat g}_1^n}+20(m-b)\varepsilon}\right)+2\exp\left(-\frac{b\ln(2)}{t_{mix}}\right).
  \end{align}
\end{prpstn}
\section{First bounds with Hoeffding-type inequality}
In this section, we will use only the boundedness property of the loss functions. We obtain bounds valid for all models, but they can be loose under particular noise conditions.
\subsection{Exponential bound}\label{sectionexponentialbound}
Let us consider out-of-sample data of length $m$ in the future of the last learning observation $X_n$. We will write a generalization bound by taking into account the last $m-b$ validation data. By doing so, we omit to take into account the $b$ first observations, but these observations account for at most $\frac{b}{m}$ in the empirical validation error because the loss function $L$ is bounded by $1$. So, we get: 
\begin{prpstn}\label{propositiongen}
  With the notations of proposition \ref{HoeffdingPaulin}, for any realization $x_1,\cdots,x_n$ of $X_1,\cdots,X_n$, integers $b$ and $m$, with $0\leq b < m$:
  \begin{equation}\label{lemgenpLbeq}
      {\mathbb P}\left(\pm\left(\hat L_{m}({\hat g}_1^n)-{\mathbb L}({\hat g}_1^n)\right)>\varepsilon+\frac{b}{m}\right)\leq 
      \exp\left(-2\frac{(m-b)\varepsilon^2}{9t_{mix}}\right)+2\exp\left(-\frac{b\ln(2)}{t_{mix}}\right).
  \end{equation}
\end{prpstn}
\subsection*{Proof} We begin to prove \eqref{lemgenpLbeq} with the sign $+$.  We have, on the event $\hat L_{m}({\hat g}_1^n)-{\mathbb L}({\hat g}_1^n)\geq\frac{b}{m}$:

\begin{align*}
  & 0\leq \frac{1}{m}\sum_{k=n+1}^{n+m}L({\hat g}_1^n(X_k))-{\mathbb L}({\hat g}_1^n)-\frac{b}{m}\leq
  \frac{1}{m}\sum_{k=n+b+1}^{n+m}L({\hat g}_1^n(X_k))-\frac{(m-b)}{m}\times {\mathbb L}({\hat g}_1^n) = \nonumber \\
  & \frac{m-b}{m}\left(\frac{1}{m-b}\sum_{k=n+b+1}^{n+m}L({\hat g}_1^n(X_k))-{\mathbb L}({\hat g}_1^n)\right)\leq
  \frac{1}{m-b}\sum_{k=n+b+1}^{n+m}L({\hat g}_1^n(X_k))-{\mathbb L}({\hat g}_1^n).
\end{align*}

So, by proposition \ref{HoeffdingPaulin}:
\begin{align*}
  & {\mathbb P}\left(\frac{1}{m}\sum_{k=n+1}^{n+m} L({\hat g}_1^n(X_k))-{\mathbb L}({\hat g}_1^n)>\varepsilon+\frac{b}{m}\right)\leq 
  {\mathbb P}\left(\frac{1}{m-b}\sum_{k=n+b+1}^{n+m}L({\hat g}_1^n(X_k))-{\mathbb L}({\hat g}_1^n)>\varepsilon\right)\leq \nonumber \\
  & \exp\left(-2\frac{(m-b)\varepsilon^2}{9t_{mix}}\right)+2\exp\left(-\frac{b\ln(2)}{t_{mix}}\right).
\end{align*}

Now, on the event $\hat L_{m}({\hat g}_1^n)-{\mathbb L}({\hat g}_1^n)<\frac{b}{m}$, we have
\[
{\mathbb P}\left(\frac{1}{m}\sum_{k=n+1}^{n+m} L({\hat g}_1^n(X_k))-{\mathbb L}({\hat g}_1^n)>\varepsilon+\frac{b}{m}\right)=0.
\]
So,
\[
{\mathbb P}\left(\hat L_{m}({\hat g}_1^n)-{\mathbb L}({\hat g}_1^n)>\varepsilon+\frac{b}{m}\right)\leq \exp\left(-2\frac{(m-b)\varepsilon^2}{9t_{mix}}\right)+2\exp\left(-\frac{b\ln(2)}{t_{mix}}\right).
\]
This shows equation \eqref{lemgenpLbeq} with the sign $+$.

For the sign $-$, remark that, on the event ${\mathbb L}({\hat g}_1^n)-\hat L_{m}({\hat g}_1^n)\geq\frac{b}{m}$:

\begin{align*}
  & 0\leq {\mathbb L}({\hat g}_1^n)-\frac{1}{m}\sum_{k=n+1}^{n+m}L({\hat g}_1^n(X_k))-\frac{b}{m}\leq
  \frac{m-b}{m}{\mathbb L}({\hat g}_1^n) -\frac{1}{m}\sum_{k=n+b+1}^{n+m}L({\hat g}_1^n(X_k)) = \nonumber \\
  & \frac{m-b}{m}\left({\mathbb L}({\hat g}_1^n) -\frac{1}{m-b}\sum_{k=n+b+1}^{n+m}L({\hat g}_1^n(X_k))\right)\leq 
  {\mathbb L}({\hat g}_1^n)-\frac{1}{m-b}\sum_{k=n+b+1}^{n+m}L({\hat g}_1^n(X_k)),
\end{align*}
  
and by the same argument as previously, we get equation \eqref{lemgenpLbeq} for the sign $-$ $\blacksquare$

Using this proposition, we can state an exponential bound for the theoretical loss function. The proof is based on a suitable choice of $b$ and is given in the Appendix.  
\begin{thrm}\label{borneL}
  With the notations of  proposition \ref{HoeffdingPaulin}, for any realization $x_1,\cdots,x_n$ of $X_1,\cdots,X_n$, and\\ $0\leq\varepsilon\leq 1$:
  \begin{equation}\label{corgeneq}
      P\left(\pm(\hat L_{m}({\hat g}_1^n)-{\mathbb L}({\hat g}_1^n))>\varepsilon\right)\leq
      \left(2\exp\left(\frac{\ln(2)}{t_{mix}}\right)+1\right)\exp\left(-\frac{m\varepsilon^2\ln(2)}{(1+9\ln(2))t_{mix}}\right).
    \end{equation}
\end{thrm}
\subsection{A first oracle inequality}\label{sectionoracle}

Let $(\left({\hat g}_1^n\right)_k)_{k=1,\cdots,N}$ denote a finite collection of prediction functions obtained by processing a realization of training sample of length $n$. In an ideal world, a benevolent oracle would tell us which index $\tilde k$ minimizes the theoretical loss:
\begin{equation}\label{theoreticalvalidationloss}
\tilde k = \arg\min_{k\in\{1,\cdots,N\}}{\mathbb L}(\left({\hat g}_1^n\right)_k).
\end{equation}
However, all we can do is to chose the index $\hat k$ that minimizes the empirical validation loss:
\begin{equation}\label{empiricalvalidationloss}
\hat k = \arg\min_{k\in\{1,\cdots,N\}}\hat L_{m}(\left({\hat g}_1^n\right)_k).
\end{equation}
An oracle inequality between the optimal and the empirical choices $\tilde k$ and $\hat k$ may be written:
\begin{equation}\label{deforacle}
{\mathbb E}\left({\mathbb L}(\left({\hat g}_1^n\right)_{\hat k})-{\mathbb L}(g^*)\right)\leq C\left({\mathbb L}(\left({\hat g}_1^n\right)_{\tilde k})-{\mathbb L}(g^*)+\frac{\gamma(n)}{n}\right).
\end{equation}
Where, $C$ is a factor at least as large as $1$, $\gamma(n)$ is a slowly growing function, and ${\mathbb L}(g^*)$ is the best expected loss.  The term $\inf_k\left({\mathbb L}\left(\left({\hat g}_1^n\right)_{k}\right)-{\mathbb L}(g^*)\right)={\mathbb L}(\left({\hat g}_1^n\right)_{\tilde k})-{\mathbb L}(g^*)$ is called the bias term, and the term $\frac{\gamma(n)}{n}$, the variance term. The concept of oracle inequality was advocated in Donoho and Johnstone \cite{Donoho} and is now widely used (see Candes \cite{Candes}).

To establish the oracle inequality,  we will begin by inequalities between empirical and theoretical losses for $\hat k$ and $\tilde k$.
The theorem \ref{borneL} and union bound give the following theorem:
\begin{thrm}\label{genboundexp}
For any realization $x_1,\cdots,x_n$ of $X_1,\cdots,X_n$, and $0\leq\varepsilon\leq 1$:
\begin{equation*}
    P\left({\mathbb L}(\left({\hat g}_1^n\right)_{\hat k})-\hat L_{m}(\left({\hat g}_1^n\right)_{\hat k})>\varepsilon\right)\leq
    N\left(2\exp\left(\frac{\ln(2)}{t_{mix}}\right)+1\right)\exp\left(-\frac{m\varepsilon^2\ln(2)}{(1+9\ln(2))t_{mix}}\right),
  \end{equation*}
and
\begin{equation*}
    P\left(\hat L_{m}\left(\left({\hat g}_1^n\right)_{\tilde k}\right)-{\mathbb L}\left(\left({\hat g}_1^n\right)_{\tilde k}\right)>\varepsilon\right)\leq
    N\left(2\exp\left(\frac{\ln(2)}{t_{mix}}\right)+1\right)\exp\left(-\frac{m\varepsilon^2\ln(2)}{(1+9\ln(2))t_{mix}}\right).
\end{equation*}
We deduce then, an  upper bound of expectations between the empirical et theoretical losses:
\begin{equation*}
    {\mathbb E}\left({\mathbb L}(\left({\hat g}_1^n\right)_{\hat k})-\hat L_{m}(\left({\hat g}_1^n\right)_{\hat k})\right) \leq \sqrt{\frac{\ln\left(eN\left(2\exp\left(\frac{\ln(2)}{t_{mix}}\right)+1\right)\right)\left(1+9\ln\left(2\right)\right)t_{mix}}{\ln\left(2\right)m}},
\end{equation*}
and
\begin{equation*}
    {\mathbb E}\left(\hat L_{m}\left(\left({\hat g}_1^n\right)_{\tilde k}\right)-{\mathbb L}\left(\left({\hat g}_1^n\right)_{\tilde k}\right))\right) \leq \sqrt{\frac{\ln\left(eN\left(2\exp\left(\frac{\ln(2)}{t_{mix}}\right)+1\right)\right)\left(1+9\ln\left(2\right)\right)t_{mix}}{\ln\left(2\right)m}}.
\end{equation*} 
\end{thrm}
\subsection*{Proof}
We can write: 
\begin{align*}
& {\mathbb E}\left({\mathbb L}(\left({\hat g}_1^n\right)_{\hat k})-\hat L_{m}(\left({\hat g}_1^n\right)_{\hat k})\right)\leq
{\mathbb E}\max\left(\left({\mathbb L}(\left({\hat g}_1^n\right)_{\hat k})-\hat L_{m}(\left({\hat g}_1^n\right)_{\hat k})\right),0\right)\leq \nonumber \\
& \sqrt{{\mathbb E}\max\left(\left({\mathbb L}(\left({\hat g}_1^n\right)_{\hat k})-\hat L_{m}(\left({\hat g}_1^n\right)_{\hat k})\right),0\right)^2}= \sqrt{\int_0^1P\left({\mathbb L}(\left({\hat g}_1^n\right)_{\hat k})-\hat L_{m}(\left({\hat g}_1^n\right)_{\hat k})>\sqrt{t}\right)dt}\leq \nonumber \\
& \sqrt{\frac{\ln\left(eN\left(2\exp\left(\frac{\ln(2)}{t_{mix}}\right)+1\right)\right)\left(1+9\ln\left(2\right)\right)t_{mix}}{\ln\left(2\right)m}}.
\end{align*}
The proof of the second inequality is symmetric.
$\blacksquare$

Now, remark that,
\begin{equation*}
    {\mathbb L}\left(\left({\hat g}_1^n\right)_{\hat k}\right)-{\mathbb L}\left(\left({\hat g}_1^n\right)_{\tilde k}\right)=
    {\mathbb L}(\left({\hat g}_1^n\right)_{\hat k})-\hat L_{m}(\left({\hat g}_1^n\right)_{\hat k})+\hat L_{m}(\left({\hat g}_1^n\right)_{\hat k})
    -\hat L_{m}\left(\left({\hat g}_1^n\right)_{\tilde k}\right)+\hat L_{m}\left(\left({\hat g}_1^n\right)_{\tilde k}\right)-{\mathbb L}\left(\left({\hat g}_1^n\right)_{\tilde k}\right),
\end{equation*}
with, by definition, $\hat L_{m}(\left({\hat g}_1^n\right)_{\hat k})-\hat L_{m}(\left({\hat g}_1^n\right)_{\tilde k})\leq 0$. Hence, we have

\begin{equation*}
      {\mathbb L}\left(\left({\hat g}_1^n\right)_{\hat k}\right)-{\mathbb L}\left(\left({\hat g}_1^n\right)_{\tilde k}\right)\leq
      \sup_{k\in\{1,\cdots,N\}}({\mathbb L}\left(\left(\hat g_1^n\right)_k\right)-\hat L_{m}\left(\left(\hat g_1^n\right)_k\right))+
      \sup_{k\in\{1,\cdots,N\}}(\hat L_{m}\left(\left(\hat g_1^n\right)_k\right)-{\mathbb L}\left(\left(\hat g_1^n\right)_k\right)),
\end{equation*}
and we get the following oracle inequality:
\begin{thrm}\label{oracle}
For any realization $(x_1,\cdots,x_n)$ of $(X_1,\cdots,X_n)$, we have:
\begin{equation*}
  {\mathbb E}\left({\mathbb L}(\left({\hat g}_1^n\right)_{\hat k})-{\mathbb L}\left(\left({\hat g}_1^n\right)_{\tilde k}\right)\right)\leq 2\sqrt{\frac{\ln\left(eN\left(2\exp\left(\frac{\ln(2)}{t_{mix}}\right)+1\right)\right)\left(1+9\ln\left(2\right)\right)t_{mix}}{\ln\left(2\right)m}}.
\end{equation*}
Or, if we denote by $g^*$ the best prediction function:
\begin{equation*}
      {\mathbb E}\left({\mathbb L}(\left({\hat g}_1^n\right)_{\hat k})-{\mathbb L}(g^*)\right)\leq {\mathbb L}\left(\left({\hat g}_1^n\right)_{\tilde k}\right)-{\mathbb L}(g^*)+2\sqrt{\frac{\ln\left(eN\left(2\exp\left(\frac{\ln(2)}{t_{mix}}\right)+1\right)\right)\left(1+9\ln\left(2\right)\right)t_{mix}}{\ln\left(2\right)m}}.
\end{equation*}
\end{thrm}
Note that  this bound rate is of the same order as in the independent case when $m$ goes toward infinite. All the bounds in this section depend on the unknown constant $t_{mix}$, however, this constant can be estimated from the data (see Wolfer and Kontorovich \cite{Wolfer}).

\section{Fast rates with Bernstein-type inequality}
We can remove the square root in the bound of theorem \ref{oracle} by increasing the empirical hold-out error estimate by a small constant factor and using Berstein-type inequality (like in Bartlett et al. \cite{Bartlett}). When the theoretical loss is small, the inequalities obtained may be better than the previous inequalities. 
\subsection{Exponential bound for ${\mathbb L}({\hat g}_1^n)$ }
We begin to establish exponential bounds for slightly modified empirical losses. As in section \ref{sectionexponentialbound}, we will write a generalization bound by considering the last $m-b$ validation data. Hence, an application of the proposition \ref{BernsteinPaulin} yields the following proposition, proven in the Appendix:
\begin{prpstn}\label{propositionlocalisation}
 With the notations and assumptions of proposition \ref{BernsteinPaulin}, let $0<a<1$ be a fixed constant. For,  any realization $(x_1,\cdots,x_n)$ of $(X_1,\cdots,X_n)$, integers $b$ and $m$, with $0\leq b < m$, we get: 
 \begin{equation}\label{Bernsteinlocal1}
     P\left(\frac{1}{1+a}\hat L_m({\hat g}_1^n)-{\mathbb L}({\hat g}_1^n)>\varepsilon+\frac{b}{m}\right)\leq
     \exp\left(-\frac{(m-b)\gamma_{ps}a(1+a)\varepsilon}{8\left(1+\frac{1}{\gamma_{ps}}\right)+20}\right)+2\exp\left(-\frac{b\ln(2)}{t_{mix}}\right),
 \end{equation}
        and
        \begin{equation}\label{Bernsteinlocal2}
	       P\left({\mathbb L}({\hat g}_1^n)-\frac{1}{1-a}\hat L_m({\hat g}_1^n)>\varepsilon+\frac{b}{m}\right)\leq
               \exp\left(-\frac{(m-b)\gamma_{ps}a(1-a)\varepsilon}{8\left(1+\frac{1}{\gamma_{ps}}\right)+20}\right)+2\exp\left(-\frac{b\ln(2)}{t_{mix}}\right).
	\end{equation}
\end{prpstn}

Using this proposition, with a suitable choice of b, we can get exponential bounds for $P\left(\frac{1}{1+a}\hat L_{m}({\hat g}_1^n)-{\mathbb L}({\hat g}_1^n)>\varepsilon\right)$ and
$P\left({\mathbb L}({\hat g}_1^n)-\frac{1}{1-a}\hat L_{m}({\hat g}_1^n)>\varepsilon\right)$. These bounds are stated in the following theorem. The proof may be found in the Appendix.  

\begin{thrm}\label{borneLbBernstein}
  With the notations and assumptions of proposition \ref{BernsteinPaulin}, for any realization $x_1,\cdots,x_n$ of $X_1,\cdots,X_n$, $0<a<1$, and $0\leq \varepsilon\leq 1$:
    \begin{equation}\label{corgeneqBernstein1}
  P\left(\frac{1}{1+a}\hat L_{m}({\hat g}_1^n)-{\mathbb L}({\hat g}_1^n)>\varepsilon\right)\leq \left(1+2\exp\left(\frac{\ln(2)}{t_{mix}}\right)\right)\exp\left(-\frac{a(1+a)m\varepsilon}{4t_{mix}\left(8\left(1+\frac{1}{\gamma_{ps}}\right)+20\right)}\right),
    \end{equation}
   and
    \begin{equation}\label{corgeneqBernstein2}
  P\left({\mathbb L}({\hat g}_1^n)-\frac{1}{1-a}\hat L_{m}({\hat g}_1^n)>\varepsilon\right)\leq \left(1+2\exp\left(\frac{\ln(2)}{t_{mix}}\right)\right)\exp\left(-\frac{a(1-a)m\varepsilon}{4t_{mix}\left(8\left(1+\frac{1}{\gamma_{ps}}\right)+20\right)}\right).
    \end{equation}
\end{thrm}

\subsection{Bound for the generalization error}\label{sectionoracleBernstein}
In the same framework as section \ref{sectionoracle}, we consider $(\left({\hat g}_1^n\right)_k)_{k=1,\cdots,N}$, a finite collection of prediction functions obtained by processing a realization of training sample of length $n$. We recall that the function with index $\hat k$ minimizes the empirical validation loss \eqref{empiricalvalidationloss}, and the function with index $\tilde k$ minimizes the theoretical loss \eqref{theoreticalvalidationloss}. 
Now, the theorem \ref{borneLbBernstein} and union bound give the following result:
\begin{thrm}\label{genboundexpBernstein}
  With the notations and assumptions of proposition \ref{BernsteinPaulin}, for any realization $x_1,\cdots,x_n$ of $X_1,\cdots,X_n$, $0<a<1$ and $0\leq \varepsilon\leq 1$:
 
\begin{equation*}
    P\left(\frac{1}{1+a}\hat L_{m}\left(\left({\hat g}_1^n\right)_{\tilde k}\right)-{\mathbb L}\left(\left({\hat g}_1^n\right)_{\tilde k}\right)>\varepsilon\right)\leq
    N\left(1+2\exp\left(\frac{\ln(2)}{t_{mix}}\right)\right)\exp\left(-\frac{a(1+a)m\varepsilon}{4t_{mix}\left(8\left(1+\frac{1}{\gamma_{ps}}\right)+20\right)}\right),
\end{equation*}
and
\begin{equation*}
    P\left({\mathbb L}(\left({\hat g}_1^n\right)_{\hat k})-\frac{1}{1-a}\hat L_{m}(\left({\hat g}_1^n\right)_{\hat k})>\varepsilon\right)\leq
    N\left(1+2\exp\left(\frac{\ln(2)}{t_{mix}}\right)\right)\exp\left(-\frac{a(1-a)m\varepsilon}{4t_{mix}\left(8\left(1+\frac{1}{\gamma_{ps}}\right)+20\right)}\right).
\end{equation*}
So we get an upper bound of the expectations of these expressions:
\begin{equation*}
    {\mathbb E}\left(\frac{1}{1+a}\hat L_{m}\left(\left({\hat g}_1^n\right)_{\tilde k}\right)-{\mathbb L}\left(\left({\hat g}_1^n\right)_{\tilde k}\right)\right) \leq \frac{4t_{mix}\left(8\left(1+\frac{1}{\gamma_{ps}}\right)+20\right)\ln\left(eN\left(2\exp\left(\frac{\ln(2)}{t_{mix}}\right)+1\right)\right)}{a(1+a)m},
\end{equation*}
and
\begin{equation*}
    {\mathbb E}\left({\mathbb L}(\left({\hat g}_1^n\right)_{\hat k})-\frac{1}{1-a}\hat L_{m}(\left({\hat g}_1^n\right)_{\hat k})\right) \leq \frac{4t_{mix}\left(8\left(1+\frac{1}{\gamma_{ps}}\right)+20\right)\ln\left(eN\left(2\exp\left(\frac{\ln(2)}{t_{mix}}\right)+1\right)\right)}{a(1-a)m}.
\end{equation*}
\end{thrm}
\subsection*{Proof} 
 We can write: 
\begin{align*}
& {\mathbb E}\left(\frac{1}{1+a}\hat L_{m}(\left({\hat g}_1^n\right)_{\tilde k}))-{\mathbb L}(\left({\hat g}_1^n\right)_{\tilde k}))\right)\leq
{\mathbb E}\max\left(\left(\frac{1}{1+a}\hat L_{m}(\left({\hat g}_1^n\right)_{\tilde k}))-{\mathbb L}(\left({\hat g}_1^n\right)_{\tilde k}))\right),0\right)= \nonumber\\
& \int_0^1P\left(\frac{1}{1+a}\hat L_{m}(\left({\hat g}_1^n\right)_{\tilde k}))-{\mathbb L}(\left({\hat g}_1^n\right)_{\tilde k}))>t\right)dt\leq\frac{4t_{mix}\left(8\left(1+\frac{1}{\gamma_{ps}}\right)+20\right)\ln\left(eN\left(2\exp\left(\frac{\ln(2)}{t_{mix}}\right)+1\right)\right)}{a(1+a)m}.
\end{align*}
The proof of the second inequality is symmetric.
$\blacksquare$

Now, remark that:
\begin{align*}
   & {\mathbb L}(\left({\hat g}_1^n\right)_{\hat k})-{\mathbb L}(\left({\hat g}_1^n\right)_{\tilde k}))=
    {\mathbb L}(\left({\hat g}_1^n\right)_{\hat k})-\frac{1}{1-a}\hat L_{m}(\left({\hat g}_1^n\right)_{\hat k})+\frac{1}{1-a}\hat L_{m}(\left({\hat g}_1^n\right)_{\hat k}) \nonumber\\
   & -\frac{1}{1-a}\hat L_{m}(\left({\hat g}_1^n\right)_{\tilde k}))+\frac{1}{1+a}\hat L_{m}(\left({\hat g}_1^n\right)_{\tilde k}))-{\mathbb L}(\left({\hat g}_1^n\right)_{\tilde k})+2\frac{a}{1-a^2}\hat L_{m}(\left({\hat g}_1^n\right)_{\tilde k},
\end{align*}
with, by definition, $\frac{1}{1-a}\hat L_{m}(\left({\hat g}_1^n\right)_{\hat k})-\frac{1}{1-a}\hat L_{m}(\left({\hat g}_1^n\right)_{\tilde k}))\leq 0$. Hence, we have

\begin{align*}
     & {\mathbb L}(\left({\hat g}_1^n\right)_{\hat k})-{\mathbb L}(\left({\hat g}_1^n\right)_{\tilde k}))\leq
      \sup_{k\in\{1,\cdots,N\}}\left({\mathbb L}\left(\left(\hat g_1^n\right)_k\right)-\frac{1}{1-a}\hat L_{m}\left(\left(\hat g_1^n\right)_k\right)\right)+\nonumber\\
    &  \sup_{k\in\{1,\cdots,N\}}\left(\frac{1}{1+a}\hat L_{m}\left(\left(\hat g_1^n\right)_k\right)-{\mathbb L}\left(\left(\hat g_1^n\right)_k\right)\right)+2\frac{a}{1-a^2}\hat L_{m}(\left({\hat g}_1^n\right)_{\tilde k},  
\end{align*}
and we get the following inequalities for the empirical choice of the index $\hat k$:
\begin{thrm}\label{oracleBernstein}
With the notations and assumptions of proposition \ref{BernsteinPaulin}, for any realization $(x_1,\cdots,x_n)$ of $(X_1,\cdots,X_n)$, and $0<a<1$, we have:
\begin{align*}
  & {\mathbb E}\left({\mathbb L}(\left({\hat g}_1^n\right)_{\hat k})-{\mathbb L}(\left({\hat g}_1^n\right)_{\tilde k}))\right)\leq\frac{4t_{mix}\left(8\left(1+\frac{1}{\gamma_{ps}}\right)+20\right)\ln\left(eN\left(2\exp\left(\frac{\ln(2)}{t_{mix}}\right)+1\right)\right)}{a(1-a)m}+\nonumber\\
 & \frac{4t_{mix}\left(8\left(1+\frac{1}{\gamma_{ps}}\right)+20\right)\ln\left(eN\left(2\exp\left(\frac{\ln(2)}{t_{mix}}\right)+1\right)\right)}{a(1+a)m}+\frac{2a}{1-a^2}{\mathbb L}(\left({\hat g}_1^n\right)_{\tilde k}).
\end{align*}
Or, if we denote by $g^*$ the best prediction function:
\begin{align*}
  & {\mathbb E}\left({\mathbb L}(\left({\hat g}_1^n\right)_{\hat k})-{\mathbb L}(g^*)\right)\leq \left(1+\frac{2a}{1-a^2}\right)\left({\mathbb L}(\left({\hat g}_1^n\right)_{\tilde k})-{\mathbb L}(g^*)\right)+\nonumber\\
  & \frac{4t_{mix}\left(8\left(1+\frac{1}{\gamma_{ps}}\right)+20\right)\ln\left(eN\left(2\exp\left(\frac{\ln(2)}{t_{mix}}\right)+1\right)\right)}{a(1-a)m}+\nonumber\\
   & \frac{4t_{mix}\left(8\left(1+\frac{1}{\gamma_{ps}}\right)+20\right)\ln\left(eN\left(2\exp\left(\frac{\ln(2)}{t_{mix}}\right)+1\right)\right)}{a(1+a)m}+\frac{2a}{1-a^2}{\mathbb L}(g^*).
\end{align*}
\end{thrm}
These bounds are not exactly like an oracle inequality \eqref{deforacle}; however, they are better than the bounds of the previous section if the theoretical loss ${\mathbb L}\left(\left({\hat g}_1^n\right)_{\tilde k}\right)$ is small enough. An extreme and rare case will be when ${\mathbb L}\left(\left({\hat g}_1^n\right)_{\tilde k}\right)=0$, and, in such case, the bounds of the previous theorem are of order $O\left(\frac{1}{m}\right)$. Note that, all the bounds in this section depend on the unknown constants $t_{mix}$ and $\gamma_{ps}$, but they can be estimated from the data (see Wolfer and Kontorovich \cite{Wolfer}).
\section{Fast rate under noise conditions}
The previous condition  ${\mathbb L}\left(\left({\hat g}_1^n\right)_{\tilde k}\right)=0$ may be seen as a rough noise condition. We can try to refine our analysis because hold-out enjoys excellent theoretical properties under noise conditions for the i.i.d. case (see Blanchard and Massart \cite{Blanchard}, Boucheron et al. \cite{Boucheron} or Massart \cite{MassartB}). Hence, we will investigate OOS properties in the Makov case under similar conditions for the noise. In this section, we assume that the state space $\cal Y$ is discrete, that $L$ is the misclassification loss function: $L(y,y')={\bf 1}_{y\neq y'}$, and the functions $g(Y_{t-1},\cdots,Y_{t-p})$ are the predictions of the following state $Y_{t}$ knowing $Y_{t-1},\cdots,Y_{t-p}$. Hence, we take advantage that for discrete observations, the prediction of the following state is a classification task. 
First, we give the assumption on the noise.
\paragraph{Assumption on the noise ${\bf (H)}$:}
\begin{itemize}
\item  A function $\omega(.)$ exists such that $\omega(x)/\sqrt{x}$ is non-increasing and, for any function $g$,
  \[
  \sqrt{Var\left({\bf 1}_{g\neq g^*}\right)}\leq \omega(\mathbb L(g)-\mathbb L(g^*)),
  \]
  where the expection are computed under the stationary law of the Markov chain $(X_t)_{t\in\mathbb Z}$.
\item Let $\tau^*_m$ denote the smallest positive solution of $\omega(\varepsilon)=\sqrt{m}\varepsilon$.
\end{itemize}
 Again, in the same framework as sections \ref{sectionoracle}, and \ref{sectionoracleBernstein}, we consider $(\left({\hat g}_1^n\right)_k)_{k=1,\cdots,N}$, a finite collection of prediction functions obtained by processing a realization of training sample of length $n$. The function with index $\hat k$ minimizes the empirical validation loss \eqref{empiricalvalidationloss}, and  the function with index $\tilde k$ minimizes the theoretical loss \eqref{theoreticalvalidationloss}.
We can then set the following proposition. We prove this proposition in the Appendix.
\begin{prpstn}\label{propositionprelNC}
  Let us assume the assumption on the noise ${\bf (H)}$, for any realization $x_1,\cdots,x_n$ of $X_1,\cdots,X_n$, integers $b$ and $m$, with $0\leq b < m$, real numbers $0\leq\varepsilon\leq 1$, and $0<\theta<1$:
 \begin{align}\label{prelNC}
    & {\mathbb P}\left({\mathbb L}\left(\left(\hat g_1^n\right)_{\hat k}\right)-{\mathbb L}(g^*)-\left(1+\theta\right)\left({\mathbb L}\left(\left(\hat g_1^n\right)_{\tilde k}\right)-{\mathbb L}(g^*)\right)>\varepsilon+\frac{(1+\theta)2b}{m}\right)\leq \nonumber\\
    & N\exp\left(-\frac{1}{1+\theta}\frac{\theta\gamma_{ps}(m-b)}{16(1+\frac{1}{\gamma_{ps}})m\tau^*_m+80\theta}\varepsilon\right)+2\exp\left(-\frac{b\ln(2)}{t_{mix}}\right).
\end{align}
\end{prpstn}
\subsection{Exponential bound under noise condition}
We can now state an exponential bound under noise condition. The following theorem is proven in the Appendix.
\begin{thrm}\label{expboundnc}
  Let us assume the assumption on the noise ${\bf (H)}$, for any realization $x_1,\cdots,x_n$ of $X_1,\cdots,X_n$, and $0<\theta<1$:
   \begin{align}\label{NC}
   & {\mathbb P}\left({\mathbb L}\left(\left(\hat g_1^n\right)_{\hat k}\right)-{\mathbb L}(g^*)-\left(1+\theta\right)\left({\mathbb L}\left(\left(\hat g_1^n\right)_{\tilde k}\right)-{\mathbb L}(g^*)\right)>\varepsilon\right)\leq\nonumber\\
   & \left(N+2\exp\left(\frac{\ln(2)}{t_{mix}}\right)\right)\exp\left(-\frac{1}{4t_{mix}(1+\theta)}\frac{\theta\gamma_{ps}m\varepsilon}{16(1+\frac{1}{\gamma_{ps}})m\tau^*_m+80\theta}\right).
   \end{align}
   Finally, taking expectation, we get the oracle inequality:
        \begin{align}\label{oracleNC}
    & {\mathbb E}\left({\mathbb L}\left(\left(\hat g_1^n\right)_{\hat k}\right)-{\mathbb L}(g^*)\right)\leq\left(1+\theta\right)\times\nonumber\\
    & \left({\mathbb L}\left(\left(\hat g_1^n\right)_{\tilde k}\right)-{\mathbb L}(g^*)+\frac{4t_{mix}(16(1+\frac{1}{\gamma_{ps}})m\tau^*_m+80\theta)}{\theta\gamma_{ps}m}\ln\left(e\left(2\exp\left(\frac{\ln(2)}{t_{mix}}\right)+N\right)\right)\right).
\end{align}
\end{thrm}
\begin{rmrk}
  Assume the Mammen-Tsybakov noise condition with exponent $\alpha$ hold (see Mammen and Tsybakov \cite{Mammen}), that is, we can choose $w(r)=\left(\frac{r}{h}\right)^{\alpha/2}$ for some positive $h$. Then, $\tau^*_m=\left(mh^\alpha\right)^{-1/(2-\alpha)}$, and the corollary translates into
     \begin{align}\label{oracleMT}
     & {\mathbb E}\left({\mathbb L}\left(\left(\hat g_1^n\right)_{\hat k}\right)-{\mathbb L}(g^*)\right)\leq\left(1+\theta\right)\times\left({\mathbb L}\left(\left(\hat g_1^n\right)_{\tilde k}\right)-{\mathbb L}(g^*)+\right.\nonumber\\
     &\left.\frac{4t_{mix}(16(1+\frac{1}{\gamma_{ps}})h^{-\alpha/(2-\alpha)}m^{1-1/(2-\alpha)}+80\theta)}{\theta\gamma_{ps}m}\ln\left(e\left(2\exp\left(\frac{\ln(2)}{t_{mix}}\right)+N\right)\right)\right)=\nonumber\\
     &\left(1+\theta\right)\times\left({\mathbb L}\left(\left(\hat g_1^n\right)_{\tilde k}\right)-{\mathbb L}(g^*)+\right.\nonumber\\
     &\left.\left(\frac{320t_{mix}}{\gamma_{ps}m}+\frac{4t_{mix}\left(16(1+\frac{1}{\gamma_{ps}})h^{-\alpha/(2-\alpha)}\right)}{\theta\gamma_{ps}m^{1/(2-\alpha)}}\right)\ln\left(e\left(2\exp\left(\frac{\ln(2)}{t_{mix}}\right)+N\right)\right)\right).
     \end{align}
\end{rmrk}
\begin{rmrk}
  If the state space $\cal Y$ is equal to $\{0,1\}$, and if the conditional expectation function:
  \[
  \eta(Y_{t-1},\cdots,Y_{t-p})={\mathbb E}\left(Y_t\left|Y_{t-1},\cdots,Y_{t-p}\right.\right)
  \]
  is such that for all $y_{t-1},\cdots,y_{t-p}\in{\cal Y}^p$, $\left|2\eta(y_{t-1},\cdots,y_{t-p})-1\right|>h$, then the  Mammen-Tsybakov noise condition holds with $\alpha=1$. The rate of the oracle inequality will then be fast:
   \begin{align*}
    & {\mathbb E}\left({\mathbb L}\left(\left(\hat g_1^n\right)_{\hat k}\right)-{\mathbb L}(g^*)\right)\leq\left(1+\theta\right)\times\left({\mathbb L}\left(\left(\hat g_1^n\right)_{\tilde k}\right)-{\mathbb L}(g^*)+\right.\\
    & \left.\left(\frac{320t_{mix}}{\gamma_{ps}m}+\frac{4t_{mix}\left(16(1+\frac{1}{\gamma_{ps}})h^{-\alpha/(2-\alpha)}\right)}{\theta\gamma_{ps}m}\right)\ln\left(e\left(2\exp\left(\frac{\ln(2)}{t_{mix}}\right)+N\right)\right)\right).
   \end{align*}  
   Hence, if ${\mathbb L}\left(\left(\hat g_1^n\right)_{\tilde k}\right)-{\mathbb L}(g^*)=0$,
     \[
     {\mathbb E}\left({\mathbb L}\left(\left(\hat g_1^n\right)_{\hat k}\right)-{\mathbb L}(g^*)\right)\leq O\left(\frac{1}{m}\right).
     \]
\end{rmrk}
\begin{rmrk}
  Let $p\geq 1$ be a fixed integer and consider the set of homogenous ergodic Markov chains of order $p$ with ${\cal Y}=\left\{0,1\right\}$.  Let us denote $\cal A$ the set of possible transition kernels:
  \begin{align*}
    {\cal A}=\left\{P(Y_{t}=1|y_{t-1},\cdots,y_{t-p}),(y_{t-1},\cdots,y_{t-p})\in\{0,1\}^p\right\}= [0,1]^{2^p}.
  \end{align*}
  Now, put the uniform measure on $\cal A$, then the Lebesgue measure of the set of models such that the  Mammen-Tsybakov noise condition does not hold with $\alpha=1$ will be null. Indeed,  $P(Y_t=0|y_{t-1},\cdots,y_{t-p})=P(Y_t=1|y_{t-1},\cdots,y_{t-p})=\frac{1}{2}$ for some $y_{t-1},\cdots,y_{t-p}$, means that a linear constraint for the coefficients of the transition kernel exists and such set is of Lebesgue measure $0$. Hence, for almost all models, an $h$ exists such that, for all $y_{t-1},\cdots,y_{t-p}\in{\cal Y}^p$, $\left|2\eta(y_{t-1},\cdots,y_{t-p})-1\right|>h$, and we get the previous rate for the oracle inequality.  
 \end{rmrk} 
\section{Appendix}
\subsection{Proof of lemma \ref{lemmemoinsdelta}}
We will prove the lemma for the sign $+$, the proof for the sign $-$ is the same. By the proposition \ref{PreBernsteinPaulin}, we have for any $0<\delta<1$:
\begin{align*}
  & P\left(E_Q(S)-S>\sqrt{\frac{8}{\gamma_{ps}}(n+1/\gamma_{ps})V_f\log\left(\frac{1}{\delta}\right)}+\frac{20}{\gamma_{ps}}B\log\left(\frac{1}{\delta}\right)\right)\leq \nonumber
  \\
  & \exp\left(-\frac{\gamma_{ps}\left(\sqrt{\frac{8}{\gamma_{ps}}(n+1/\gamma_{ps})V_f\log\left(\frac{1}{\delta}\right)}+\frac{20}{\gamma_{ps}}B\log\left(\frac{1}{\delta}\right)\right)^2}{8\left(n+1/\gamma_{ps}\right)V_f+20B\left(\sqrt{\frac{8}{\gamma_{ps}}(n+1/\gamma_{ps})V_f\log\left(\frac{1}{\delta}\right)}+\frac{20}{\gamma_{ps}}B\log\left(\frac{1}{\delta}\right)\right)}\right)\leq \nonumber
  \\
  & \exp\left(-\frac{8(n+1/\gamma_{ps})V_f\log\left(\frac{1}{\delta}\right)+40\sqrt{\frac{8}{\gamma_{ps}}(n+1/\gamma_{ps})V_f\log\left(\frac{1}{\delta}\right)}B\log\left(\frac{1}{\delta}\right)+\frac{1}{\gamma_{ps}}\left(20B\log\left(\frac{1}{\delta}\right)\right)^2}{8\left(n+1/\gamma_{ps}\right)V_f+20B\left(\sqrt{\frac{8}{\gamma_{ps}}(n+1/\gamma_{ps})V_f\log\left(\frac{1}{\delta}\right)}+\frac{20}{\gamma_{ps}}B\log\left(\frac{1}{\delta}\right)\right)}\right)\leq \nonumber
  \\
 & \exp\left(-\log\left(\frac{1}{\delta}\right)\frac{8(n+1/\gamma_{ps})V_f+40B\sqrt{\frac{8}{\gamma_{ps}}(n+1/\gamma_{ps})V_f\log\left(\frac{1}{\delta}\right)}+\frac{1}{\gamma_{ps}}\left(20B\right)^2\log\left(\frac{1}{\delta}\right)}{8\left(n+1/\gamma_{ps}\right)V_f+20B\left(\sqrt{\frac{8}{\gamma_{ps}}(n+1/\gamma_{ps})V_f\log\left(\frac{1}{\delta}\right)}+\frac{20}{\gamma_{ps}}B\log\left(\frac{1}{\delta}\right)\right)}\right)\leq \delta.
\end{align*}

Noting that $\frac{8}{\gamma_{ps}}(n+1/\gamma_{ps})=\frac{8}{\gamma^2_{ps}}(\gamma_{ps}n+1)\leq \frac{8(\gamma_{ps}+1)}{\gamma^2_{ps}}n$, completes the proof $\blacksquare$

\subsection{Proof of theorem \ref{borneL}}
We prove \eqref{corgeneq} with the sign $+$, the proof for the sign $-$ is symmetric.
Since $0\leq \varepsilon\leq 1$, if $b =\lfloor \frac{m\varepsilon^2}{1+9\ln(2)} \rfloor$, then $\varepsilon-\frac{b}{m}>0$. Using the proposition \ref{propositiongen}, we have:
\[
  {\mathbb P}\left(\pm(\hat L_{m}({\hat g}_1^n)-{\mathbb L}({\hat g}_1^n))>\varepsilon\right)\leq \exp\left(-2\frac{(m-b)(\varepsilon-\frac{b}{m})^2}{9t_{mix}}\right)+2\exp\left(-\frac{b\ln(2)}{t_{mix}}\right).
 \]
 Now, if ${\tilde b} = \frac{m\varepsilon^2}{1+9\ln(2)}$, then ${\tilde b}-1\leq b\leq {\tilde b}$, and we get
 \[
   {\mathbb P}\left(\pm(\hat L_{m}({\hat g}_1^n)-{\mathbb L}({\hat g}_1^n))>\varepsilon\right)\leq
   \exp\left(-2\frac{(m-{\tilde b})(\varepsilon-\frac{\tilde b}{m})^2}{9t_{mix}}\right)+2\exp\left(\frac{\ln(2)}{t_{mix}}\right)\exp\left(-\frac{{\tilde b}\ln(2)}{t_{mix}}\right).
 \]
  Moreover,
  \begin{align*}
   & \exp\left(-2\frac{(m-{\tilde b})(\varepsilon-\frac{\tilde b}{m})^2}{9t_{mix}}\right)+2\exp\left(\frac{\ln(2)}{t_{mix}}\right)\exp\left(-\frac{{\tilde b}\ln(2)}{t_{mix}}\right)=\\
   & \left(2\exp\left(\frac{\ln(2)}{t_{mix}}\right)+\exp\left(-m\left(\frac{2(1-\frac{\tilde b}{m})(\varepsilon-\frac{\tilde b}{m})^2}{9t_{mix}}-\frac{\tilde b\ln(2)}{mt_{mix}}\right)\right)\right)\exp\left(-\frac{\tilde b\ln(2)}{t_{mix}}\right)=\\
   & \left(2\exp\left(\frac{\ln(2)}{t_{mix}}\right)+\exp\left(-\frac{m}{t_{mix}}\left(\frac{2(1-\frac{\tilde b}{m})(\varepsilon-\frac{\tilde b}{m})^2}{9}-\frac{\tilde b\ln(2)}{m}\right)\right)\right)\exp\left(-\frac{\tilde b\ln(2)}{t_{mix}}\right).
  \end{align*}
  Now,
  \begin{align*}
  &\frac{2(1-\frac{\tilde b}{m})(\varepsilon-\frac{\tilde b}{m})^2}{9}-\frac{\tilde b\ln(2)}{m}=
    \frac{2(1-\frac{\varepsilon^2}{1+9\ln(2)})(\varepsilon-\frac{\varepsilon^2}{1+9\ln(2)})^2}{9}-\frac{\varepsilon^2\ln(2)}{1+9\ln(2)}\geq\\
  &  \varepsilon^2\left(\frac{2(1-\frac{1}{1+9\ln(2)})(1-\frac{1}{1+9\ln(2)})^2}{9}-\frac{\ln(2)}{1+9\ln(2)}\right)=\\
  &  \varepsilon^2\left(\frac{2}{9}\left(1-\frac{1}{1+9\ln(2)}\right)^3-\frac{\ln(2)}{1+9\ln(2)}\right)>0,
    \end{align*}
  and finally
  \begin{align*}
    &\exp\left(-2\frac{(m-{\tilde b})(\varepsilon-\frac{\tilde b}{m})^2}{9t_{mix}}\right)+2\exp\left(\frac{\ln(2)}{t_{mix}}\right)\exp\left(-\frac{{\tilde b}\ln(2)}{t_{mix}}\right)\leq\\
    &\left(2\exp\left(\frac{\ln(2)}{t_{mix}}\right)+1\right)\exp\left(-\frac{\tilde b\ln(2)}{t_{mix}}\right)=\\
    &\left(2\exp\left(\frac{\ln(2)}{t_{mix}}\right)+1\right)\exp\left(-\frac{m\varepsilon^2\ln(2)}{(1+9\ln(2))t_{mix}}\right).
  \end{align*}
 $\blacksquare$

\subsection{Proof of proposition \ref{propositionlocalisation}}
First, we will prove that for a stationary uniformly ergodic Markov chain $(X_t)_{t\in \mathbb Z}$ and any function $g\in \cal G$:
\begin{equation}\label{PrelimBernstein}
  P\left(\frac{1}{1+a}\frac{1}{m}\sum_{k=1}^m L(g(X_k))-{\mathbb E}(L(g(X)))>\varepsilon\right)\leq \exp\left(- \frac{m\gamma_{ps}a(1+a)\varepsilon}{8\left(1+\frac{1}{\gamma_{ps}}\right)+20}\right),
\end{equation}
where ${\mathbb E}(L(g(X)))$ is computed under the stationay law of $(X_t)_{t\in \mathbb Z}$.

Since $0\leq L(g(X))\leq 1$,
$0\leq L(g(X))^2 \leq L(g(X))\leq 1$, and
$V(L(g(X)))\leq {\mathbb E}(L(g(X)))(1-{\mathbb E}(L(g(X))))\leq {\mathbb E}(L(g(X)))$.

Moreover,
\begin{align*}
  &P\left(\frac{1}{1+a}\frac{1}{m}\sum_{k=1}^m L(g(X_k))-{\mathbb E}(L(g(X)))>\varepsilon\right)=\\
  &P\left(\frac{1}{m}\sum_{k=1}^m L(g(X_k))-{\mathbb E}(L(g(X)))>a{\mathbb E}(L(g(X)))+(1+a)\varepsilon\right).
\end{align*}

Let $t=a {\mathbb E}(L(g(X)))+(1+a)\varepsilon$, by the proposition \ref{PreBernsteinPaulin}, we have:
\begin{align*}
      &P\left(\frac{1}{m}\sum_{k=1}^{m}L(g(X_k))>t\right)\leq \exp\left(-\frac{m^2t^2\gamma_{ps}}{8(m+1/\gamma_{ps}){\mathbb E}(L(g(X)))+20 m t}\right)=\\
      &\exp\left(-\frac{mt^2\gamma_{ps}}{8(1+\frac{1}{m\gamma_{ps}}){\mathbb E}(L(g(X)))+20 t}\right)
\end{align*}

Now,
\[
8(1+\frac{1}{m\gamma_{ps}}){\mathbb E}(L(g(X)))+20  t\leq t\left(\frac{8(1+1/\gamma_{ps})}{a}+20\right),
\]
hence
\begin{align*}
      &P\left(\frac{1}{m}\sum_{k=1}^{m}L(g(X_k))>t\right)\leq \exp\left(-\frac{mt\gamma_{ps}}{\frac{8(1+1/\gamma_{ps})}{a}+20}\right)\leq\\
      &\exp\left(-\frac{mat\gamma_{ps}}{8(1+1/\gamma_{ps})+20}\right)\leq\exp\left(-\frac{ma(1+a)\varepsilon\gamma_{ps}}{8(1+1/\gamma_{ps})+20}\right),
\end{align*}
and we deduce equation \eqref{PrelimBernstein}.

Now, using  equation (3.27) of Paulin \cite{Paulin}, we get
  \begin{equation}\label{FastBernstein}
      P\left(\frac{1}{1+a}\frac{1}{m-b}\sum_{k=n+b+1}^{n+m}L({\hat g}_1^n(X_k))-{\mathbb L}({\hat g}_1^n)>\varepsilon\right)\leq\nonumber
      \exp\left(-\frac{(m-b)\gamma_{ps}a(1+a)}{8\left(1+\frac{1}{\gamma_{ps}}\right)+20}\varepsilon\right)+2\exp\left(-\frac{b\ln(2)}{t_{mix}}\right).
  \end{equation}
  Finally, by the same arguments as in the proof of proposition \ref{propositiongen}, we get equation \eqref{Bernsteinlocal1}. The proof of equation \eqref{Bernsteinlocal2} is symmetric.
$\blacksquare$

\subsection{Proof of theorem \ref{borneLbBernstein}}
We prove \eqref{corgeneqBernstein1}, the proof for \eqref{corgeneqBernstein2} is symmetric.
Let us define ${\tilde b}=\frac{ma(1+a)\varepsilon}{4\ln(2)\left(8\left(1+\frac{1}{\gamma_{ps}}\right)+20\right)}$. Since $\lfloor {\tilde b}\rfloor> {\tilde b}-1$, equation \eqref{Bernsteinlocal1} yields
\begin{align*}
  &P\left(\frac{1}{1+a}\hat L_{m}({\hat g}_1^n)-{\mathbb L}({\hat g}_1^n)>\varepsilon\right)\leq\\
  &\exp\left(-\frac{\gamma_{ps}a(1+a)}{8\left(1+\frac{1}{\gamma_{ps}}\right)+20}\left(m-\frac{\tilde b}{m}\right)\left(\varepsilon-\frac{\tilde b}{m}\right)\right)+2\exp\left(\frac{\ln(2)}{t_{mix}}\right)\exp\left(-\frac{{\tilde b}\ln(2)}{t_{mix}}\right)=\\
  &\exp\left(-\frac{\gamma_{ps}a(1+a)m\varepsilon}{8\left(1+\frac{1}{\gamma_{ps}}\right)+20}\left(1-\frac{a(1+a)}{4\ln(2)\left(8\left(1+\frac{1}{\gamma_{ps}}\right)+20\right)}\varepsilon\right)\left(1-\frac{a(1+a)}{4\ln(2)\left(8\left(1+\frac{1}{\gamma_{ps}}\right)+20\right)}\right)\right)\\
  &+2\exp\left(\frac{\ln(2)}{t_{mix}}\right)\exp\left(-\frac{a(1+a)m\varepsilon}{4t_{mix}\left(8\left(1+\frac{1}{\gamma_{ps}}\right)+20\right)}\right)=\\
  &\left.\exp\left(-\frac{\gamma_{ps}a(1+a)m\varepsilon}{8\left(1+\frac{1}{\gamma_{ps}}\right)+20}\left(\left(1-\frac{a(1+a)}{4\ln(2)\left(8\left(1+\frac{1}{\gamma_{ps}}\right)+20\right)}\varepsilon\right)\left(1-\frac{a(1+a)}{4\ln(2)\left(8\left(1+\frac{1}{\gamma_{ps}}\right)+20\right)}\right)-\frac{1}{4t_{mix}\gamma_{ps}}\right)\right)\right)\times\\
  &\exp\left(-\frac{a(1+a)m\varepsilon}{4t_{mix}\left(8\left(1+\frac{1}{\gamma_{ps}}\right)+20\right)}\right).
\end{align*}

Now, since $\varepsilon\leq 1$,  and, by equation (3.9) of Paulin \cite{Paulin}, $\gamma_{ps}\geq \frac{1}{2t_{mix}}\Leftrightarrow 2\geq\frac{1}{2t_{mix}\gamma_{ps}}$, we get
\begin{align*}
  &P\left(\frac{1}{1+a}\hat L_{m}({\hat g}_1^n)-{\mathbb L}({\hat g}_1^n)>\varepsilon\right)\leq\\
  &\left(2\exp\left(\frac{\ln(2)}{t_{mix}}\right)+\exp\left(-\gamma_{ps}\frac{a(1+a)m\varepsilon}{8\left(1+\frac{1}{\gamma_{ps}}\right)+20}\left(\left(1-\frac{a(1+a)}{4\ln(2)\left(8\left(1+\frac{1}{\gamma_{ps}}\right)+20\right)}\right)^2-\frac{1}{2}\right)\right)\right)\times\\
  &\exp\left(-\frac{a(1+a)m\varepsilon}{4t_{mix}\left(8\left(1+\frac{1}{\gamma_{ps}}\right)+20\right)}\right).
\end{align*}
Finally, noting that $\frac{a(1+a)}{8\left(1+\frac{1}{\gamma_{ps}}\right)+20}\leq \frac{1}{14}$, we have
\(
\left(1-\frac{a(1+a)}{4\ln(2)\left(8\left(1+\frac{1}{\gamma_{ps}}\right)+20\right)}\right)^2-\frac{1}{2}\geq 0,
\)
and
\[
  P\left(\frac{1}{1+a}\hat L_{m}({\hat g}_1^n)-{\mathbb L}({\hat g}_1^n)>\varepsilon\right)\leq \left(2\exp\left(\frac{\ln(2)}{t_{mix}}\right)+1\right)\exp\left(-\frac{a(1+a)m\varepsilon}{4t_{mix}\left(8\left(1+\frac{1}{\gamma_{ps}}\right)+20\right)}\right)
\]
 $\blacksquare$ 
\subsection{Proof of proposition \ref{propositionprelNC}}
Let us considere a finite collection of functions $\{g_1,\cdots,g_N\}$. For any integers $b$ and $m$, with $0\leq b<m$, and a sample $(X_1,\cdots,X_{m+b})$ where the $m$ last variables  $(X_{b+1},\cdots,X_{m+b})$ follow the stationary law of the Markov chain $(X_t)_{t\in \mathbb Z}$. Let us define:
\[
g_{\hat k}=\min_{k\in\{1,\cdots,N\}} \frac{1}{m+b}\sum_{t=1}^{m+b}L(g_k(X_t))\mbox{ and }g_{\tilde k}=\min_{k\in\{1,\cdots,N\}}{\mathbb L}(g_k).
\]
We will give bounds involving the $m$ last variables: ${\hat L}_m(g_k):=\frac{1}{m}\sum_{t=b+1}^{m+b}L(g_k(X_t))$.
Note that, for any function $g_k$:
\begin{equation}\label{relation_m_mb1}
    \frac{1}{m+b}\sum_{t=1}^{m+b}L(g_k(X_t))-\frac{1}{m}\sum_{t=b+1}^{m+b}L(g_k(X_t))\leq 
    \frac{1}{m+b}\sum_{t=1}^{m+b}L(g_k(X_t))-\frac{1}{m+b}\sum_{t=b+1}^{m+b}L(g_k(X_t))\leq \frac{b}{m+b},
\end{equation}
and
\begin{align}\label{relation_m_mb2}
    &\frac{1}{m}\sum_{t=b+1}^{m+b}L(g_k(X_t))-\frac{1}{m+b}\sum_{t=1}^{m+b}L(g_k(X_t))\leq 
    \frac{m+b}{m(m+b)}\sum_{t=b+1}^{m+b}L(g_k(X_t))-\frac{m}{m(m+b)}\sum_{t=1}^{m+b}L(g_k(X_t))\leq\nonumber\\
    &\frac{b}{m(m+b)}\sum_{t=b+1}^{m+b}L(g_k(X_t))\leq \frac{b}{m+b},
\end{align}
so
\begin{equation}\label{relation_hatk}
    \frac{1}{m}\sum_{t=b+1}^{m+b}L(g_{\hat k}(X_t))-\frac{1}{m}\sum_{t=b+1}^{m+b}L(g_k(X_t))\leq \frac{2b}{m+b}.
\end{equation}

By the lemma \ref{lemmemoinsdelta} and the union bound, with probability at least $1-\delta$, for all $k\in \{1,\cdots,N\}$,
\[
  {\mathbb L}\left(g_{k}\right)-{\mathbb L}(g^*)\leq
  {\hat L}_m\left(g_{k}\right)-{\hat L}_m(g^*)+\sqrt{\frac{8(1+\frac{1}{\gamma_{ps}})\log\left(\frac{N}{\delta}\right)}{\gamma_{ps}m}}\times\omega\left({\mathbb L}(g_k)-{\mathbb L}(g^*)\right)+\frac{40\log\left(\frac{N}{\delta}\right)}{\gamma_{ps}m},
\]
and
\[
  {\mathbb L}(g^*)-{\mathbb L}\left(g_{\tilde k}\right)\leq
  {\hat L}_m(g^*)-{\hat L}_m\left(g_{\tilde k}\right)+\sqrt{\frac{8(1+\frac{1}{\gamma_{ps}})\log\left(\frac{N}{\delta}\right)}{\gamma_{ps}m}}\times\omega\left({\mathbb L}\left(g_{\tilde k}\right)-{\mathbb L}(g^*)\right)+\frac{40\log\left(\frac{N}{\delta}\right)}{\gamma_{ps}m}.
\]
Since ${\mathbb L}\left(g_{\tilde k}\right)-{\mathbb L}(g^*)\leq {\mathbb L}\left(g_{k}\right)-{\mathbb L}(g^*)$, for any $k\in \{1,\cdots,N\}$, by summing the two inequalities, we obtain
\begin{equation*}
  {\mathbb L}\left(g_{k}\right)-{\mathbb L}\left(g_{\tilde k}\right)\leq{\hat L}_m\left(g_{k}\right)-{\hat L}_m\left(g_{\tilde k}\right)+
  2\sqrt{\frac{8(1+\frac{1}{\gamma_{ps}})\log\left(\frac{N}{\delta}\right)}{\gamma_{ps}m}}\times\omega\left({\mathbb L}\left(g_{\tilde k}\right)-{\mathbb L}(g^*)\right)+\frac{80\log\left(\frac{N}{\delta}\right)}{\gamma_{ps}m}.
\end{equation*}
As ${\hat L}_m\left(g_{\hat k}\right)-{\hat L}_m\left(g_{\tilde k}\right)\leq \frac{2b}{m+b}$, with probability larger than $1-\delta$,
\begin{equation*}
    {\mathbb L}\left(g_{\hat k}\right)-{\mathbb L}\left(g_{\tilde k}\right)\leq
    \frac{2b}{m+b}+2\sqrt{\frac{8(1+\frac{1}{\gamma_{ps}})\log\left(\frac{N}{\delta}\right)}{\gamma_{ps}m}}\times\omega\left({\mathbb L}\left(g_{\tilde k}\right)-{\mathbb L}(g^*)\right)+\frac{80\log\left(\frac{N}{\delta}\right)}{\gamma_{ps}m}.
\end{equation*}
Let $\tau^*_m$ be defined as the statement of the theorem. If ${\mathbb L}\left(g_{\hat k}\right)-{\mathbb L}(g^*)\geq \tau^*_{m+b}$, then $\omega\left({\mathbb L}\left(g_{\hat k}\right)-{\mathbb L}(g^*)\right)/\sqrt{m+b}\leq \sqrt{{\mathbb L}\left(g_{\hat k}\right)-{\mathbb L}(g^*)\tau^*_{m+b}}$, and we have
\begin{equation*}
    {\mathbb L}\left(g_{\hat k}\right)-{\mathbb L}\left(g_{\tilde k}\right)\leq
     \frac{2b}{m+b}+2\sqrt{\frac{8(1+\frac{1}{\gamma_{ps}})\log\left(\frac{N}{\delta}\right)}{\gamma_{ps}}}\times \sqrt{\frac{m+b}{m}}\sqrt{\tau^*_{m+b}}\sqrt{{\mathbb L}\left(g_{\hat k}\right)-{\mathbb L}(g^*)}+\frac{80\log\left(\frac{N}{\delta}\right)}{\gamma_{ps}m}.
\end{equation*}
For $0<\theta<1$ we have:
\begin{align*}
  &\frac{\theta^2}{2}\left({\mathbb L}\left( g_{\hat k}\right)-{\mathbb L}(g^*)\right)
   -2\sqrt{\frac{8(1+\frac{1}{\gamma_{ps}})\log\left(\frac{N}{\delta}\right)}{\gamma_{ps}}}\times \sqrt{\frac{m+b}{m}}\sqrt{\tau^*_{m+b}}\sqrt{{\mathbb L}\left(g_{\hat k}\right)-{\mathbb L}(g^*)}\theta\\
  & +\frac{16(1+\frac{1}{\gamma_{ps}})\log\left(\frac{N}{\delta}\right)}{\gamma_{ps}}\frac{m+b}{m}\tau^*_{m+b}=
  \left(\frac{\theta}{\sqrt{2}}\sqrt{{\mathbb L}\left(g_{\hat k}\right)-{\mathbb L}(g^*)}-\sqrt{\frac{16(1+\frac{1}{\gamma_{ps}})\log\left(\frac{N}{\delta}\right)}{\gamma_{ps}}}\times \sqrt{\frac{m+b}{m}}\sqrt{\tau^*_{m+b}}\right)^2\geq 0,
\end{align*}
and
\begin{equation*}
    2\sqrt{\frac{8(1+\frac{1}{\gamma_{ps}})\log\left(\frac{N}{\delta}\right)}{\gamma_{ps}}}\times \sqrt{\frac{m+b}{m}}\sqrt{\tau^*_{m+b}}\sqrt{{\mathbb L}\left(g_{\hat k}\right)-{\mathbb L}(g^*)}\leq
    \frac{\theta}{2}\left({\mathbb L}\left(g_{\hat k}\right)-{\mathbb L}(g^*)\right)+\frac{16(1+\frac{1}{\gamma_{ps}})\log\left(\frac{N}{\delta}\right)}{\theta\gamma_{ps}}\frac{m+b}{m}\tau^*_{m+b},
\end{equation*}
so
\begin{equation*}
    {\mathbb L}\left(g_{\hat k}\right)-{\mathbb L}\left(g_{\tilde k}\right)\leq \frac{\theta}{2}\left({\mathbb L}\left(g_{\hat k}\right)-{\mathbb L}(g^*)\right)
+\frac{2b}{m+b}+\frac{16(1+\frac{1}{\gamma_{ps}})\log\left(\frac{N}{\delta}\right)}{\theta\gamma_{ps}}\frac{m+b}{m}\tau^*_{m+b}+\frac{80\log\left(\frac{N}{\delta}\right)}{\gamma_{ps}m}.
\end{equation*}
Hence, with probability larger than $1-\delta$
\begin{equation*}
    \left(1-\frac{\theta}{2}\right)\left({\mathbb L}\left(g_{\hat k}\right)-{\mathbb L}(g^*)\right)\leq{\mathbb L}\left(g_{\tilde k}\right)-{\mathbb L}(g^*)
+\frac{2b}{m+b}+\frac{16(1+\frac{1}{\gamma_{ps}})\log\left(\frac{N}{\delta}\right)}{\theta\gamma_{ps}}\frac{m+b}{m}\tau^*_{m+b}+\frac{80\log\left(\frac{N}{\delta}\right)}{\gamma_{ps}m},
\end{equation*}
and
\begin{equation*}
    \left({\mathbb L}\left(g_{\hat k}\right)-{\mathbb L}(g^*)\right)\leq
  \frac{1}{1-\frac{\theta}{2}}\left({\mathbb L}\left(g_{\tilde k}\right)-{\mathbb L}(g^*)+\frac{2b}{m+b}+\frac{16(1+\frac{1}{\gamma_{ps}})\log\left(\frac{N}{\delta}\right)}{\theta\gamma_{ps}}\frac{m+b}{m}\tau^*_{m+b}+\frac{80\log\left(\frac{N}{\delta}\right)}{\gamma_{ps}m}\right).
\end{equation*}
Since, for $0<\theta<1$, $\frac{1}{1-\frac{\theta}{2}}\leq 1+\theta$, we get, with probability larger than $1-\delta$:
\begin{equation*}
    {\mathbb L}\left(g_{\hat k}\right)-{\mathbb L}(g^*)\leq\left(1+\theta\right)\times
    \left({\mathbb L}\left(g_{\tilde k}\right)-{\mathbb L}(g^*)+\frac{2b}{m+b}+\frac{16(1+\frac{1}{\gamma_{ps}})\log\left(\frac{N}{\delta}\right)}{\theta\gamma_{ps}}\frac{m+b}{m}\tau^*_{m+b}+\frac{80\log\left(\frac{N}{\delta}\right)}{\gamma_{ps}m}\right),
\end{equation*}
or
\begin{equation}\label{borneenprobadelta}
    {\mathbb L}\left(g_{\hat k}\right)-{\mathbb L}(g^*)-\left(1+\theta\right)\left({\mathbb L}\left(g_{\tilde k}\right)-{\mathbb L}(g^*)+\frac{2b}{m+b}\right)\leq
  \left(1+\theta\right)\left(\frac{16(1+\frac{1}{\gamma_{ps}})(m+b)\tau^*_{m+b}+80\theta}{\theta\gamma_{ps}m}\log\left(\frac{N}{\delta}\right)\right).
\end{equation}
Note that we have done the reasoning if ${\mathbb L}\left(g_{\hat k}\right)-{\mathbb L}(g^*)\geq \tau^*_{m+b}$. However, if ${\mathbb L}\left(g_{\hat k}\right)-{\mathbb L}(g^*)< \tau^*_{m+b}$, the bound \eqref{borneenprobadelta} is obvious.
Now, we deduce from equation \eqref{borneenprobadelta} that
\begin{equation*}
    {\mathbb P}\left({\mathbb L}\left(g_{\hat k}\right)-{\mathbb L}(g^*)-\left(1+\theta\right)\left({\mathbb L}\left(g_{\tilde k}\right)-{\mathbb L}(g^*)+\frac{2b}{m+b}\right)>\varepsilon\right)\leq
    N\exp\left(-\frac{1}{1+\theta}\frac{\theta\gamma_{ps}m\varepsilon}{16(1+\frac{1}{\gamma_{ps}})(m+b)\tau^*_{m+b}+80\theta}\right).
  \end{equation*}
 Now, considering the actual chain $(X_t)_{t\in\mathbb N}$, in the framework of section \ref{sectionoracle},  we get for any realization $x_1,\cdots,x_n$ of $X_1,\cdots,X_n$, integers $b$ and $m$, with $0\leq b < m$:
 \begin{align*}
    &{\mathbb P}\left({\mathbb L}\left(\left(\hat g_1^n\right)_{\hat k}\right)-{\mathbb L}(g^*)-\left(1+\theta\right)\left({\mathbb L}\left(\left(\hat g_1^n\right)_{\tilde k}\right)-{\mathbb L}(g^*)\right)>\varepsilon+\frac{(1+\theta)2b}{m}\right)\leq\nonumber\\
    &N\exp\left(-\frac{1}{1+\theta}\frac{\theta\gamma_{ps}(m-b)\varepsilon}{16(1+\frac{1}{\gamma_{ps}})m\tau^*_m+80\theta}\right)+2\exp\left(-\frac{b\ln(2)}{t_{mix}}\right).
\end{align*}
 $\blacksquare$
 \subsection{Proof of theorem \ref{expboundnc}}
Applying the proposition \ref{propositionprelNC} to $2(1+\theta)\varepsilon$, we get
 \begin{align*}
    &{\mathbb P}\left(\frac{{\mathbb L}\left(\left(\hat g_1^n\right)_{\hat k}\right)-{\mathbb L}(g^*)}{2(1+\theta)}-\frac{\left({\mathbb L}\left(\left(\hat g_1^n\right)_{\tilde k}\right)-{\mathbb L}(g^*)\right)}{2}>\varepsilon+\frac{b}{m}\right)\leq\nonumber\\
    &N\exp\left(-2\frac{\theta\gamma_{ps}(m-b)}{16(1+\frac{1}{\gamma_{ps}})m\tau^*_m+80\theta}\varepsilon\right)+2\exp\left(-\frac{b\ln(2)}{t_{mix}}\right).
\end{align*}
 Let us define ${\tilde b}=\frac{\theta m\varepsilon}{(16(1+\frac{1}{\gamma_{ps}})m\tau^*_m+80\theta)2\ln(2)}$.
 Following the same reasoning as in the proof of theorem \ref{borneLbBernstein}, we get
\begin{align*}
  &{\mathbb P}\left({\mathbb L}\left(\left(\hat g_1^n\right)_{\hat k}\right)-{\mathbb L}(g^*)-\left(1+\theta\right)\left({\mathbb L}\left(\left(\hat g_1^n\right)_{\tilde k}\right)-{\mathbb L}(g^*)\right)>\varepsilon\right)\leq\\
  &\left(2\exp\left(\frac{\ln(2)}{t_{mix}}\right)+N\exp\left(-\frac{2\theta\gamma_{ps}m\varepsilon}{(16(1+\frac{1}{\gamma_{ps}})m\tau^*_m+80\theta)}\left(\left(1-\frac{\theta}{(16(1+\frac{1}{\gamma_{ps}})m\tau^*_m+80\theta)2\ln(2)}\right)^2-\frac{1}{2}\right)\right)\right)\times\\
  &\exp\left(-\frac{\theta\gamma_{ps}\varepsilon m}{(16(1+\frac{1}{\gamma_{ps}})m\tau^*_m+80\theta)2t_{mix}}\right).
\end{align*}
Noting that $\frac{\theta}{16(1+\frac{1}{\gamma_{ps}})m\tau^*_m+80\theta}\leq \frac{1}{80}$, we have
\[
\left(1-\frac{\theta}{(16(1+\frac{1}{\gamma_{ps}})m\tau^*_m+80\theta)2\ln(2)}\right)^2-\frac{1}{2}\geq 0,
\]
and
\begin{align*}
  &{\mathbb P}\left({\mathbb L}\left(\left(\hat g_1^n\right)_{\hat k}\right)-{\mathbb L}(g^*)-\left(1+\theta\right)\left({\mathbb L}\left(\left(\hat g_1^n\right)_{\tilde k}\right)-{\mathbb L}(g^*)\right)>\varepsilon\right)\leq\\
  &\left(2\exp\left(\frac{\ln(2)}{t_{mix}}\right)+N\right)\exp\left(-\frac{1}{4t_{mix}(1+\theta)}\frac{\theta\gamma_{ps}\varepsilon m}{16(1+\frac{1}{\gamma_{ps}})m\tau^*_m+80\theta}\right).
\end{align*}
$\blacksquare$
\bibliographystyle{plain}
\bibliography{biblio.bib}

\end{document}